\documentclass{article}
\usepackage{amsmath,amsthm,amsopn,amstext,amscd,amsfonts,amssymb}
\usepackage{natbib}

\def\diag{\mathop{\rm diag}\nolimits}

\def\build#1#2#3{\mathrel{\mathop{#1}\limits^{#2}_{#3}}}
\def\e{\mathop{\rm e}\nolimits}
\def\tr{\mathop{\rm tr}\nolimits}

\def\vecp{\mathop{\rm vecp}\nolimits}

\def\etr{\mathop{\rm etr}\nolimits}

\renewenvironment{abstract}
                 {\vspace{6pt}
                  \begin{center}
                  \begin{minipage}{5in}
                  \centerline{\textbf{Abstract}}
                  \noindent\ignorespaces
                 }
                 {\end{minipage}\end{center}}
\newtheorem{thm}{\textbf{Theorem}}[section]
\newtheorem{cor}{\textbf{Corollary}}[section]
\newtheorem{lem}{\textbf{Lemma}}[section]
\theoremstyle{definition}

\title{\Large \textbf{Shape Theory Via SV Decomposition II}}
\author{
  \textbf{Jos\'e A. D\'{\i}az-Garc\'{\i}a} \thanks{Corresponding author\newline
    }\\
  {\normalsize Department of Statistics and Computation} \\
  {\normalsize Universidad Aut\'onoma Agraria Antonio Narro}\\
  {\normalsize 25350 Buenavista, Saltillo, Coahuila, Mexico} \\
  {\normalsize E-mail: jadiaz@uaaan.mx} \\[2ex]
  \textbf{Francisco J. Caro-Lopera} \\
  {\normalsize Department of Basic Sciences} \\
  {\normalsize Universidad de Medell\'{\i}n} \\
  {\normalsize Carrera 87 No.30-65, of. 5-103}\\
  {\normalsize Medell\'{\i}n, Colombia}\\
  {\normalsize E-mail: fjcaro@udem.edu.co}\\
}
\date{}
\begin{document}
\maketitle

\begin{abstract}
The non isotropic and  non central elliptical shape distributions via the Le and
Kendall SVD decomposition approach are derived in this paper in the context of
invariant polynomials and zonal polynomials. The so termed cone and disk densities
here obtained generalise some results of the literature. Finally, some particular
densities are applied in a classical data of Biology, and the inference is performed
after choosing the best model by using a modified BIC criterion.
\end{abstract}

\section{Introduction and the main principle.}

The multivariate statistical shape theory has been developed in the last two decades
around the classical works based on normality and isotropy, see \citet{GM93} and
\citet{DM98} and the references there in. Recent works extended this results to
elliptical models and partial non isotropy, see for example \citet{Caro2009} and
\citet{dgr:03}, however some important problems remain, the study of the shape theory
without any restriction of the covariance matrix in the elliptical model.

The problem arises from the point of view of applications, the isotropic assumption
$\boldsymbol{\Theta}=\mathbf{I}_{K}$ for an elliptical shape model of the form
$$
  \mathbf{X} \sim \mathcal{E}_{N \times K} (\boldsymbol{\mu}_{{}_{\mathbf{X}}},
    \boldsymbol{\Sigma}_{{}_{\mathbf{X}}},\boldsymbol{\Theta}, h),
$$
restricts substantially the correlations of the landmarks in the figure.  Then, we
expect the non isotropic model, with any positive definite matrix
$\boldsymbol{\Theta}$, as the best model for considering all the possible
correlations among the anatomical (geometrical o mathematical) points.

This problem can be solved by considering the following procedure: Let be
$$
  \mathbf{X} \sim \mathcal{E}_{N \times K} (\boldsymbol{\mu}_{{}_{\mathbf{X}}},
    \boldsymbol{\Sigma}_{{}_{\mathbf{X}}},\boldsymbol{\Theta}, h),
$$
if $\boldsymbol{\Theta}^{1/2}$ is the positive definite square root of the matrix
$\boldsymbol{\Theta}$, i .e. $\boldsymbol{\Theta} = (\boldsymbol{\Theta}^{1/2})^{2}$,
with $\boldsymbol{\Theta}^{1/2}:$ $K \times K$, \citet[p. 11]{gv:93}, and noting that
$$
  \mathbf{X} \boldsymbol{\Theta}^{-1} \mathbf{X}' = \mathbf{X}
  (\boldsymbol{\Theta}^{-1/2}\boldsymbol{\Theta}^{-1/2})^{-1}\mathbf{X}' = \mathbf{X}
  \boldsymbol{\Theta}^{-1/2} (\mathbf{X} \boldsymbol{\Theta}^{-1/2})' =
  \mathbf{Z}\mathbf{Z}',
$$
where
$$
  \mathbf{Z} = \mathbf{X} \boldsymbol{\Theta}^{-1/2},
$$
then
$$
  \mathbf{Z} \sim \mathcal{E}_{N \times K}(\boldsymbol{\mu}_{{}_{\mathbf{Z}}},
  \boldsymbol{\Sigma}_{{}_{\mathbf{X}}}, \mathbf{I}_{K},h)
$$
with $\boldsymbol{\mu}_{{}_{\mathbf{Z}}} = \boldsymbol{\mu}_{{}_{\mathbf{X}}}
\boldsymbol{\Theta}^{-1/2}$, (see \citet[p. 20]{gv:93}).

And we arrive at the classical starting point in shape theory where the original
landmark matrix is replaced by $\mathbf{Z} = \mathbf{X} \boldsymbol{\Theta}^{-1/2}$.
Then we can proceed as usual, removing from $\mathbf{Z}$, translation, scale,
rotation and/or reflection in order to obtain the shape of $\mathbf{Z}$ (or
$\mathbf{X}$) via the QR, SVD, polar decompositions. The shape theory associated with
the SVD decomposition can be study from two different approaches, one due to
\citet{g:91} and another proposed by \citet{LK93}.

We study in this paper the statistical approach of Le and Kendall under a generalised
elliptical model.  First, recall some facts of this technique (\citet{LK93}). It is
known that the shape of an object is all geometrical information that remains after
filtering out translation, rotation and scale information of an original figure
(represented by a matrix $\mathbf{X}$) comprised in $N$ landmarks in $K$ dimensions.
Hence, we say that two figures, $\mathbf{X}_{1}:N\times K$ and
$\mathbf{X}_{2}:N\times K$ have the same shape if they are related with a special
similarity transformation $\mathbf{X}_{2}=\beta
\mathbf{X}_{1}\mathbf{H}+\mathbf{1}_{N}\boldsymbol{\gamma}'$, where
$\mathbf{H}:K\times K\in SO(K)$ (the rotation), $\boldsymbol{\gamma}:K\times 1$ (the
translation), $\mathbf{1}_{N}:N\times 1,$ $\mathbf{1}_{N}=(1,1,\ldots,1)'$, and
$\beta>0$ (the scale). Thus, in this context, the shape of a matrix $\mathbf{X}$ is
all the geometrical information about $\mathbf{X}$ that is invariant under Euclidean
similarity transformations. Then, the shape space is the set of all possible shapes,
it is the orbit space of the non-coincident $N$ landmarks in $\Re^{K}$ under the
action of the Euclidean similarity transformations. The dimension of this space is
$NK-K-1-K(K-1)/2$, it is, the original dimension $NK$ is reduced by $K$ for location,
by 1 for uniform scale and by $K(K-1)/2$ for rotation. In other words, the shape of
$\mathbf{X}$ is the set $\{\mathbf{P}\boldsymbol{\Gamma}:\boldsymbol{\Gamma}\in
SO(K)\}$ where $\mathbf{P}$ is the so termed pre-shape of $\mathbf{X}$ defined as
$\mathbf{P}=\mathbf{L}\mathbf{X}/\|\mathbf{L}\mathbf{X}\|$ ($\mathbf{L}$ is Helmert
submatrix, for example) which is invariant under translation and scaling of
$\mathbf{X}$. The rotated $\mathbf{P}$ on the pre-shape sphere  is termed a fibre of
the pre-shape space $S_{K}^{N}$, these fibres do not overlap and corresponds one to
one with shapes in the shape space $\Sigma_{K}^{N}$, it is, the pre-shape space is
partitioned into fibres by the rotation group $SO(K)$ and the fibre is the orbit of
$\mathbf{P}$ under the action of $SO(K)$. Thus, $\Sigma_{K}^{N}$ is the quotient
space of $\Sigma_{K}^{N}$ under the action of $SO(K)$, in notation
$\Sigma_{K}^{N}=S_{K}^{N}/SO(K)$, which means that the shape of $\mathbf{X}$ is an
equivalent class under the action of the group of similarity transformations.

Now, the statistical theory of shape associated to this approach studies the effect
of randomness and assume a probabilistic model for the original matrix in order to
obtain the density of the pre-shape (cone) and shape (disk). The complete procedure
for obtaining the shape of an original $\mathbf{X}$ can be summarised in the
following steps:
$$
  \mathbf{L}\mathbf{X}\boldsymbol{\Theta}^{-1/2}=\mathbf{L}\mathbf{Z} =
  \mathbf{Y}=\mathbf{V}'\mathbf{D}\mathbf{H}= r \mathbf{V}'\mathbf{W}
  \mathbf{H}=r \mathbf{V}'\mathbf{W}(\mathbf{u})\mathbf{H},
$$
where the matrix $\mathbf{L}:(N-1)\times N$ has orthonormal rows to
$\textbf{1}=(1,\ldots,1)'$.  $\mathbf{L}$ can be a submatrix of the Helmert matrix,
for example. Here $\mathbf{Y}=\mathbf{V}'\mathbf{D}\mathbf{H}$  is the SVD of matrix
$\mathbf{Y}$, with $\mathbf{V}:n\times (N-1)$ and $\mathbf{H}:n\times K$
semiorthogonal matrices and $\mathbf{D}:n\times n,$
$\mathbf{D}=\diag(D_{1},\ldots,D_{n})$; $\mathbf{W}=\mathbf{D}/r$,
$r=\|\mathbf{D}\|=\left(\sum_{i=1}^{n}D_{i}^{2}\right)^{1/2}=\|\mathbf{Y}\|$.

Then the standard problem considers a model for $\mathbf{X}$ and finds the so termed
cone and disk densities, which are the densities of $\mathbf{D}$ and
$\mathbf{W}(\mathbf{u})$, respectively.  A number of published works have studied the
classical problem which assumes a non central Gaussian model for $X$ very restricted
by the isotropy assumption $\mathbf{\Theta} = \mathbf{I}_{K}$ and
$\mathbf{\Sigma}_{\mathbf{X}} = \sigma^{2} \mathbf{I}_{N}$ in
$$
   \mathbf{X} \sim \mathcal{N}_{N \times K} (\boldsymbol{\mu}_{{}_{\mathbf{X}}},
   \boldsymbol{\Sigma}_{{}_{\mathbf{X}}},\boldsymbol{\Theta}).
$$
These restrictions facilitates the integration in terms of the zonal polynomials and
the asymptotic densities can be derived and applied. This procedure based on
normality and isotropy is very common in literature of shape under QR, SVD and affine
decompositions too, see \citet{GM93} and \citet{DM98}, and the references there in.

However, it is clear that the gaussian case do not support all the applications and
the statistical theory of shape could enriched if complete families of cone and disk
densities are available for a particular experiment and the researcher can model the
situation by applying a model selection criteria (see for example \citet{ri:78},
\citet{kr:95}, \citet{r:95} and \citet{YY07}, among many others).

Therefore, in this paper we propose the statistical theory of shape of Le and
Kendall's approach under any non central elliptical model without any restriction of
the covariance matrix.

Explicitly, consider a full  covariance elliptical model, indexed by
the generator function $h$,
$$
   \mathbf{X} \sim \mathcal{E}_{N \times K} (\boldsymbol{\mu}_{{}_{\mathbf{X}}},
   \boldsymbol{\Sigma}_{{}_{\mathbf{X}}},\boldsymbol{\Theta}, h),
$$
so, by the main principle given  above, at the beginning of this
section, we have that
$$
  \mathbf{Z} \sim \mathcal{E}_{N \times K}(\boldsymbol{\mu}_{{}_{\mathbf{Z}}},
  \boldsymbol{\Sigma}_{{}_{\mathbf{X}}}, \mathbf{I}_{K},h),
$$
with $\mathbf{Z} = \mathbf{X} \boldsymbol{\Theta}^{-1/2}$,
$\boldsymbol{\mu}_{{}_{\mathbf{Z}}} = \boldsymbol{\mu}_{{}_{\mathbf{X}}}
\boldsymbol{\Theta}^{-1/2}$.

Then the (SVD) Le and Kendall's  shape coordinates $\mathbf{u}$ of $\mathbf{X}$ are
constructed in several steps summarised in the expression
\begin{equation}\label{eq:SVDLEKSteps}
  \mathbf{L}\mathbf{X}\boldsymbol{\Theta}^{-1/2}=\mathbf{L}\mathbf{Z} =
  \mathbf{Y}=\mathbf{V}'\mathbf{D}\mathbf{H}= r \mathbf{V}'\mathbf{W}\mathbf{H}=r
  \mathbf{V}'\mathbf{W}(\mathbf{u})\mathbf{H}.
\end{equation}

Denote $\boldsymbol{\mu}=\mathbf{L}\boldsymbol{\mu}_{{}_{\mathbf{X}}}$, so
$\mathbf{Y}:(N-1)\times K$ is invariant to translations of the figure $\mathbf{Z}$,
and
$$
  \mathbf{Y}\sim\mathcal{E}_{N-1\times K}(\boldsymbol{\mu} \boldsymbol{\Theta}^{-1/2},
    \boldsymbol{\Sigma}\otimes \mathbf{I}_{K},h),
$$
where $\boldsymbol{\Sigma}=\mathbf{L}\boldsymbol{\Sigma}_{\mathbf{X}}\mathbf{L}'$,
meanwhile the matrix $\mathbf{W}$, the shape of $\mathbf{X}$, is invariant under
(translations), rotations and scaling of the landmark data matrix $\mathbf{X}$.

By considering the above main principe,  this paper  studies  Le and Kendall's
approach for the shape theory based on the SVD decomposition and any non isotropic
non central elliptical model. Section \ref{sec:SVDII} obtains the general densities,
the so termed cone and disk densities, with some corollaries. Then, the central case
and its invariance is studied in Section \ref{sec:Centralconedisk}. At the end
inference on small and large mouse vertebra data is performed with the classical
Gaussian model and two non normal Kotz models, then the best model is chosen by a
modified BIC criterion and the corresponding test for equality in mean disk shape is
obtained.

\section{Shape Theory via SVD  Le and Kendall's approach} \label{sec:SVDII}

We start with the   jacobian of the corresponding decomposition.

\begin{lem}\label{lem:SVDjacobian}
Let be $\mathbf{Y}:N-1\times K$, then there exist $\mathbf{V}\in \mathcal{V}_{n,N-1}$
represents the Stiefel manifold, $\mathbf{H}\in \mathcal{V}_{n,K}$ and
$\mathbf{D}:n\times n$, $\mathbf{D}=\diag (D_{1},\ldots,D_{n})$, $n=\min [(N-1),K]$;
$D_{1}\geq D_{2}\geq \cdots \geq D_{n}\geq 0$, such that
$\mathbf{Y}=\mathbf{V}'\mathbf{D}\mathbf{H}$; This factorisation is termed
non-singular part of the SVD. Then
\begin{equation*}
    (d\mathbf{Y})=2^{-n}|\mathbf{D}|^{N-1+K-2n}\displaystyle\prod_{i<j}^{n}(D_{i}^{2}-D_{j}^{2})
    (d\mathbf{D})(\mathbf{V}d\mathbf{V}')(\mathbf{H}d\mathbf{H}').
\end{equation*}
\end{lem}
\textit{Proof.} See \citet{dgm:97}. \qed

In order to obtain the joint density function of $(V,D)$ we need the following
generalisation of \citet[eq. (22)]{JAT64}.

\begin{lem}\label{lem:geneq22James64}
Let $\mathbf{X}: K \times n$, $\mathbf{Y}: K \times K$ and $\mathbf{H} \in
\mathcal{V}_{n,K}$. Then
\begin{enumerate}
  \item
  $$
    \int_{\mathbf{H} \in \mathcal{V}_{n,K}} [\tr(\mathbf{Y} + \mathbf{X}
    \mathbf{H})]^{p} (\mathbf{H}d\mathbf{H}') = \displaystyle\frac{2^{n}
    \pi^{Kn/2}}{\Gamma_{n}[\frac{1}{2} K]} \sum_{f=0}^{\infty} \sum_{\lambda}
    \displaystyle\frac{(p)_{2f} (\tr \mathbf{Y})^{p-2f}}{(\frac{1}{2} K)_{\lambda}}
    \displaystyle\frac{C_{\lambda}(\frac{1}{4} \mathbf{X}\mathbf{X}')}{f!},
  $$
  where $|(\tr \mathbf{Y})^{-1} \tr \mathbf{X}\mathbf{H}|<1$ and $\tr \mathbf{Y} \neq 0$.
  \item
  $
    \displaystyle \int_{\mathbf{H} \in \mathcal{V}_{n,K}} \tr(\mathbf{Y} + \mathbf{X}\mathbf{H})
    \etr(r(\mathbf{Y} + \mathbf{X}\mathbf{H}))(\mathbf{H}d\mathbf{H}') =
  $
  \par \noindent \hfill
    \hbox{$\displaystyle\frac{2^{n}
    \pi^{Kn/2}}{\Gamma_{n}[\frac{1}{2} K]} \etr(r \mathbf{Y}) \left \{ \tr \mathbf{Y}
    {}_{0}F_{1} (\frac{1}{2} K; \displaystyle\frac {r^{2}}{4} \mathbf{X}\mathbf{X}') + \sum_{f=0}^{\infty}
    \sum_{\lambda} \displaystyle\frac{(f+\frac{1}{2})}{(\frac{1}{2} K)_{\lambda}}
    \displaystyle\frac{C_{\lambda}(\frac{1}{4} \mathbf{X}\mathbf{X}')}{f!} \right \}, $}
  \par \noindent
\end{enumerate}
where $p \in \Re$, $r \in \Re$, $C_{\kappa}(\mathbf{B})$ are the zonal polynomials of
$\mathbf{B}$ corresponding to the partition $\kappa=(f_{1},\ldots f_{p})$ of $f$,
with $\sum_{i=1}^{p}f_{i}=f$; and $(a)_{\kappa}=\prod_{i=1}(a-(j-1)/2)_{f_{j}}$,
$(a)_{f}=a(a+1)\cdots (a+f-1)$, are the generalised hypergeometric coefficients and
${}_{0}F_{1}$ is the Bessel function, \citet{JAT64}.
\end{lem}

\textit{Proof.}
\begin{enumerate}
   \item From Lemma 9.5.3 \citet[Lemma 9.5.3, p. 397]{MR1982} we have
      $$
        \int_{\mathbf{H} \in \mathcal{V}_{n,K}} [\tr(\mathbf{Y} + \mathbf{X}
        \mathbf{H})]^{p} (\mathbf{H}d\mathbf{H}') = \frac{2^{n}
        \pi^{Kn/2}}{\Gamma_{n}[\frac{1}{2} K]} \int_{\mathcal{O}(K)}[\tr(\mathbf{Y} +
        \mathbf{X}\mathbf{H})]^{p} (d\mathbf{H}).
      $$
      Furthermore, for $\tr \mathbf{Y} \neq 0$ and $|(\tr \mathbf{Y})^{-1} \tr \mathbf{X}\mathbf{H}| <1$
      $$
        [\tr(\mathbf{Y} + \mathbf{X}\mathbf{H})]^{p} = (\tr \mathbf{Y})^{p} \sum_{f=0}^{\infty}
        \frac{(p)_{f}}{f!}(\tr \mathbf{Y})^{-f}(\tr \mathbf{X}\mathbf{H})^{f}.
      $$
      Now from \citet[eqs. (46) and (22)]{JAT64}) it follows that\\[2ex]
      $
        \int_{\mathbf{H} \in \mathcal{V}_{n,K}} [\tr(\mathbf{Y} + \mathbf{X}\mathbf{H})]^{p}
        (\mathbf{H}d\mathbf{H}') =
      $
      \par \noindent \hfill
    \hbox{$\displaystyle\frac{2^{n}
        \pi^{Kn/2}}{\Gamma_{n}[\frac{1}{2} K]} \sum_{f=0}^{\infty}
        \sum_{\lambda} \frac{(p)_{2f} (\tr \mathbf{Y})^{-2f}}{(2f)!}
        \frac{(\frac{1}{2})_{f}}{(\frac{1}{2} K)_{\lambda}} C_{\lambda}(\mathbf{X}\mathbf{X}'),$}
    \par \noindent\newline

      the result follows, noting that $(\frac{1}{2})_{f}/(2f)! = 1/(4^{f}
      f!)$ and that $C_{\lambda}(a\mathbf{X}\mathbf{X}') = a^{f} C_{\lambda}(\mathbf{X}\mathbf{X}')$.
   \item This follows by expanding the exponentials in series
      of powers and by applying (22) and (27) from \citet{JAT64}.
\end{enumerate}
\qed

Thus, we can obtain:

\begin{thm}\label{th:jointVD}
The joint density of $(\mathbf{V},\mathbf{D})$ is
\begin{eqnarray*}
    f_{\mathbf{V},\mathbf{D}}(\mathbf{V},\mathbf{D})&=&\frac{\pi^{\frac{nK}{2}}
    |\mathbf{D}|^{N-1+K-2n}\displaystyle\prod_{i<j}(D_{i}^{2}-D_{j}^{2})}{\Gamma_{n}\left[\frac{1}{2}
    K\right]|\boldsymbol{\Sigma}|^{\frac{K}{2}}}\\
    &&\times\sum_{t=0}^{\infty}\sum_{\kappa}\frac{h^{(2t)}\left[\tr
    \left(\boldsymbol{\Sigma}^{-1} \mathbf{V}'\mathbf{D}^{2}\mathbf{V} +
    \boldsymbol{\Omega}\right)\right] } {t!\left(\frac{1}{2}K\right)_{\kappa}}C_{\kappa}
    \left(\boldsymbol{\Omega}\boldsymbol{\Sigma}^{-1}\mathbf{V}'\mathbf{D}^{2}\mathbf{V}\right).
\end{eqnarray*}
\end{thm}
\textit{Proof.} Let be $\boldsymbol{\Omega} = \boldsymbol{\Sigma}^{-1}
\boldsymbol{\mu}\boldsymbol{\Theta}^{-1} \boldsymbol{\mu}',$ then the density of
$\mathbf{Y}$ is given by
$$
  f_{\mathbf{Y}}(\mathbf{Y})=\frac{1}{|\boldsymbol{\Sigma}|^{\frac{K}{2}}} h\left[\tr
    \left(\boldsymbol{\Sigma}^{-1}\mathbf{Y}\mathbf{Y}' + \boldsymbol{\Omega}\right) - 2
    \tr \boldsymbol{\mu}'\boldsymbol{\Sigma}^{-1}\mathbf{Y}\right].
$$
Now, make the change of variables $\mathbf{Y}=\mathbf{V}'\mathbf{D}\mathbf{H}$, so, by Lemma
\ref{lem:SVDjacobian}, the joint density function of $\mathbf{V}$, $\mathbf{D}$, $\mathbf{H}$
is
\begin{eqnarray*}
    dF_{\mathbf{V},\mathbf{D},\mathbf{H}}(\mathbf{V},\mathbf{D},\mathbf{H}) &=&
    \frac{2^{-n} |\mathbf{D}|^{N-1+K-2n} \displaystyle
    \prod_{i<j}(D_{i}^{2}-D_{j}^{2})}{|\boldsymbol{\Sigma}|^{\frac{K}{2}}}
    (\mathbf{V}d\mathbf{V}')(d\mathbf{D})\\
    && \times h\left[\tr\left(\boldsymbol{\Sigma}^{-1}\mathbf{V}'\mathbf{D}^{2}\mathbf{V}
    + \boldsymbol{\Omega}\right)-2\tr\boldsymbol{\mu}' \boldsymbol{\Sigma}^{-1}
    \mathbf{V}'\mathbf{D}\mathbf{H}\right](\mathbf{H}d\mathbf{H}').
\end{eqnarray*}
Expanding in power series
\begin{eqnarray*}
    && dF_{\mathbf{V},\mathbf{D},\mathbf{H}}(\mathbf{V},\mathbf{D},\mathbf{H}) =
    \frac{2^{-n}|\mathbf{D}|^{N-1+K-2n}\displaystyle
    \prod_{i<j}(D_{i}^{2}-D_{j}^{2})}{|\boldsymbol{\Sigma}|^{\frac{K}{2}}}
    (\mathbf{V}d\mathbf{V}')(d\mathbf{D})\\
    &&\quad\times \sum_{t=0}^{\infty}\frac{1}{t!} h^{(t)}\left[\tr
    \left(\boldsymbol{\Sigma}^{-1}\mathbf{V}'
    \mathbf{D}^{2}\mathbf{V}+\boldsymbol{\Omega}\right)\right]
    \left[\tr\left(-2\tr\boldsymbol{\mu}'\boldsymbol{\Sigma}^{-1}\mathbf{V}'
    \mathbf{D}\mathbf{H}\right)\right]^{t}(\mathbf{H}d\mathbf{H}').
\end{eqnarray*}
From Lemma \ref{lem:geneq22James64}
\begin{footnotesize}
\begin{equation*}
    \int_{\mathbf{H} \in\mathcal{V}_{n,K}}\left[\tr\left(-2\tr\boldsymbol{\mu}'
    \boldsymbol{\Sigma}^{-1}\mathbf{V}'\mathbf{D}\mathbf{H}\right)\right]^{2t}
    (\mathbf{H}d\mathbf{H}') = \frac{2^{n}\pi^{\frac{nK}{2}}}{\Gamma_{n}\left[\frac{1}{2}
    K\right]} \sum_{\kappa} \frac{\left(\frac{1}{2}\right)_{t}
    4^{t}}{\left(\frac{1}{2}K\right)_{\kappa}}
    C_{\kappa}\left(\boldsymbol{\Omega}\boldsymbol{\Sigma}^{-1}\mathbf{V}'
    \mathbf{D}^{2}\mathbf{V}\right).
\end{equation*}
\end{footnotesize}
Observing that $\frac{\left(\frac{1}{2}\right)_{t}4^{t}}{(2t)!}=\frac{1}{t!}$, the
marginal joint density of $\mathbf{V}$, $\mathbf{D}$ is given by
\begin{eqnarray*}
    && dF_{\mathbf{V},\mathbf{D}}(\mathbf{V},\mathbf{D}) = \frac{\pi^{\frac{nK}{2}}
    |\mathbf{D}|^{N-1+K-2n} \displaystyle \prod_{i<j}(D_{i}^{2}-D_{j}^{2})}
    {\Gamma_{n}\left(\frac{K}{2}\right)|\boldsymbol{\Sigma}|^{\frac{K}{2}}}\\
    &&\quad\times \sum_{t=0}^{\infty}\frac{h^{(2t)}
    \left[\tr\left(\boldsymbol{\Sigma}^{-1}\mathbf{V}' \mathbf{D}^{2}\mathbf{V} +
    \boldsymbol{\Omega}\right)\right]}{t!\left(\frac{1}{2}K\right)_{\kappa}}
    C_{\kappa}\left(\boldsymbol{\Omega}\boldsymbol{\Sigma}^{-1}\mathbf{V}'
    \mathbf{D}^{2}\mathbf{V}\right)(\mathbf{V}d\mathbf{V}')(d\mathbf{D}). \qed
\end{eqnarray*}

Now, note that $\mathbf{D}$ contains $n$ coordinates for which, under this method,
its corresponding joint density is termed cone density (or size-and-shape density).
Then we have the first main result of this section.

\begin{thm}\label{th:conedensity}
The cone density is given by
\begin{eqnarray}\label{eq:conedensity}
    f_{\mathbf{D}}(\mathbf{D}) &=& \frac{2^{n}\pi^{\frac{n(N-1+K)}{2}}
    |\mathbf{D}|^{N-1+K-2n}\displaystyle\prod_{i<j}\left(D_{i}^{2}-D_{j}^{2}\right)}
    {\Gamma_{n}\left[\frac{1}{2} K\right]\Gamma_{n}\left[\frac{1}{2}
    (N-1)\right]|\boldsymbol{\Sigma}|^{\frac{K}{2}}}\nonumber\\
    && \times\sum_{\theta,\kappa}\sum_{\phi\in\theta\cdot\kappa}
    \frac{h^{(2t+l)}(\tr\boldsymbol{\Omega})\Delta_{\phi}^{\theta,\kappa}
    C_{\phi}\left(\mathbf{D}^{2}\right)C_{\phi}^{\theta,\kappa}
    \left(\boldsymbol{\Sigma}^{-1},\boldsymbol{\Omega}\boldsymbol{\Sigma}^{-1}\right)}
    {t!l!\left(\frac{1}{2}K\right)_{\kappa}C_{\phi}\left(\mathbf{I}_{N-1}\right)},
\end{eqnarray}
where the notation of the sum operators, $C_{\phi}^{\theta,\kappa}$ and
$\Delta_{\phi}^{\theta,\kappa}$ are given in \citet{d:80}, in particular
$\Delta_{\phi}^{\theta,\kappa}=
\displaystyle\frac{C_{\phi}^{\theta,\kappa}(\mathbf{I},\mathbf{I})}{C_{\phi}(\mathbf{I})}$.
\end{thm}
\textit{Proof.} The joint density of $\mathbf{V}$, $\mathbf{D}$ is
\begin{eqnarray*}
    dF_{\mathbf{V},\mathbf{D}}(\mathbf{V},\mathbf{D}) &=& \frac{\pi^{\frac{nK}{2}}
    |\mathbf{D}|^{N-1+K-2n}\displaystyle\prod_{i<j}\left(D_{i}^{2}-D_{j}^{2}\right)}
    {\Gamma_{n} \left[\frac{1}{2}K\right]|\boldsymbol{\Sigma}|^{\frac{K}{2}}}
    (d\mathbf{D}) (\mathbf{V}d\mathbf{V}')\\
    && \times\sum_{t=0}^{\infty}\sum_{\kappa}\frac{h^{(2t)} \left[\tr
    \left(\boldsymbol{\Sigma}^{-1}\mathbf{V}' \mathbf{D}^{2} \mathbf{V} +
    \boldsymbol{\Omega}\right)\right]C_{\kappa}\left(\boldsymbol{\Omega}
    \boldsymbol{\Sigma}^{-1}\mathbf{V}'\mathbf{D}^{2}
    \mathbf{V}\right)}{t!\left(\frac{1}{2}K\right)_{\kappa}}.
\end{eqnarray*}
Assuming that $h^{(2t)}(\cdot)$ can be expanded in power series,
\begin{eqnarray*}
    h^{(2t)} \left[\tr\boldsymbol{\Sigma}^{-1}\mathbf{V}'\mathbf{D}^{2}
    \mathbf{V}+\tr\boldsymbol{\Omega}\right] &=& \sum_{l=0}^{\infty} \frac{h^{(2t+l)}
    \left[\tr\boldsymbol{\Omega}\right]}{l!} \left[\tr\boldsymbol{\Sigma}^{-1}\mathbf{V}'
    \mathbf{D}^{2} \mathbf{V}\right]^{l}\\
    &=&\sum_{l=0}^{\infty}\sum_{\theta}\frac{h^{(2t+l)}\left[\tr
    \boldsymbol{\Omega}\right]}{l!}C_{\theta}\left(\boldsymbol{\Sigma}^{-1}\mathbf{V}'
    \mathbf{D}^{2}\mathbf{V}\right),
\end{eqnarray*}
where $C_{\theta}(\mathbf{A})$ is the zonal polynomial corresponding to the partition
$\theta=(l_{1},\ldots,l_{\alpha})$, with $\sum_{i=1}^{\alpha}l_{i}=l$.

From \citet{d:80}, eq. (4.13), the integration of $dF_{\mathbf{V},
\mathbf{D}}(\mathbf{V}, \mathbf{D})$ with respect to $\mathbf{V}\in
\mathcal{V}_{n,N-1}$ results
\begin{eqnarray*}
    && \int_{\mathbf{V}\in \mathcal{V}_{n,N-1}}C_{\kappa}\left(\boldsymbol{\Omega}
    \boldsymbol{\Sigma}^{-1}\mathbf{V}'\mathbf{D}^{2}\mathbf{V}\right)
    C_{\theta}\left(\boldsymbol{\Sigma}^{-1}\mathbf{V}'\mathbf{D}^{2}\mathbf{V}\right)
    (\mathbf{V}d\mathbf{V}')\\
    && \quad\quad\quad =\frac{2^{n}\pi^{\frac{n(N-1)}{2}}}{\Gamma_{n}\left[\frac{1}{2}
    (N-1)\right]} \sum_{\phi\in\theta\cdot\kappa}
    \frac{C_{\phi}^{\theta,\kappa}\left(\boldsymbol{\Sigma}^{-1},
    \boldsymbol{\Omega}\boldsymbol{\Sigma}^{-1}\right)
    C_{\phi}^{\theta,\kappa}\left(\mathbf{D}^{2},\mathbf{D}^{2}\right)}
    {C_{\phi}\left(\mathbf{I}_{N-1}\right)}.
\end{eqnarray*}
And by \citet{d:80}, eq. (5.1) we have that
\begin{eqnarray*}
    f_{\mathbf{D}}(\mathbf{D}) &=& \frac{2^{n}\pi^{\frac{n(N-1+K)}{2}}
    |\mathbf{D}|^{N-1+K-2n} \displaystyle\prod_{i<j}\left(D_{i}^{2}-D_{j}^{2}\right)}
    {\Gamma_{n}\left[\frac{1}{2} K\right]\Gamma_{n}\left[\frac{1}{2} (N-1)\right]
    |\boldsymbol{\Sigma}|^{\frac{K}{2}}}\\
    &=&\sum_{\theta,\kappa}\sum_{\phi\in\theta\cdot\kappa} \frac{h^{(2t)}
    \left(\tr\boldsymbol{\Omega}\right)\Delta_{\phi}^{\theta,\kappa}C_{\phi}\left(\mathbf{D}^{2}\right)
    C_{\phi}^{\theta,\kappa}\left(\boldsymbol{\Sigma}^{-1}, \boldsymbol{\Omega}
    \boldsymbol{\Sigma}^{-1}\right)} {t!l!\left(\frac{1}{2}K\right)_{\kappa}
    C_{\phi}\left(\mathbf{I}_{N-1}\right)}. \qed
\end{eqnarray*}

Now, we can derive the isotropic version of the cone density

\begin{cor}\label{cor:diskisotropic}
Let be $\boldsymbol{\Sigma}=\sigma^{2}\mathbf{I}$, then
\begin{eqnarray}\label{eq:coneIsotropic}
    f_{\mathbf{D}}(\mathbf{D})&=&\frac{2^{n}\pi^{\frac{n(N-1+K)}{2}} \displaystyle
    \prod_{i=1}^{n} D_{i}^{N-1+K-2n} \prod_{i<j}\left(D_{i}^{2}-D_{j}^{2}\right)}
    {\Gamma_{n}\left[\frac{1}{2} K\right]\Gamma_{n}\left[\frac{1}{2}
    (N-1)\right]\sigma^{(N-1)K}}\nonumber\\
    & & \times\sum_{t=0}^{\infty}\sum_{\kappa} \frac{h^{(2t)}(\tr\boldsymbol{\Omega } +
    \frac{1}{\sigma^{2}} \tr \mathbf{D}^{2})C_{\kappa} \left(\frac{1}{\sigma^{2}}
    \mathbf{D}^{2}\right) C_{\kappa} \left(\boldsymbol{\Omega} \right)}{t!
    \left(\frac{1}{2}K\right)_{\kappa} C_{\kappa}\left(\mathbf{I}\right)}.
\end{eqnarray}
\end{cor}
\textit{Proof.}  From Theorem \ref{th:conedensity}
\begin{enumerate}
    \item  $\boldsymbol{\Omega}=\boldsymbol{\Sigma}^{-1}\boldsymbol{\mu}
    \boldsymbol{\Theta}^{-1}\boldsymbol{\mu}'= \displaystyle\frac{1}{\sigma^{2}}\boldsymbol{\mu}
    \boldsymbol{\Theta}^{-1}\boldsymbol{\mu}'$.
    \item From \citet{d:80}, eq. (5.7),
    \begin{eqnarray*}
        C_{\phi}^{\theta,\kappa}\left(\boldsymbol{\Sigma}^{-1},\boldsymbol{\Omega}
        \boldsymbol{\Sigma}^{-1}\right)&=& C_{\phi}^{\theta,\kappa}\left(\frac{1}{\sigma^{2}}
        \mathbf{I}_{N-1},\frac{1}{\sigma^{2}}\boldsymbol{\Omega}\right),\\
        &=&\left(\frac{1}{\sigma^{2}}\right)^{2t+l}C_{\phi}^{\theta,\kappa}\left(\mathbf{I}_{N-1},
        \boldsymbol{\Omega}\right),\\
        &=&\left(\frac{1}{\sigma^{2}}\right)^{2t+l}\frac{\Delta_{\phi}^{\theta,\kappa}
        C_{\phi}(\mathbf{I}_{N-1}) C_{\kappa}(\boldsymbol{\Omega})}{C_{\kappa}(\mathbf{I}_{N-1})}.
\end{eqnarray*}
\end{enumerate}
Therefore the second line of (\ref{eq:conedensity}), denoted by $J$, it is simplified
as follows:
\begin{eqnarray*}
    J &=& \sum_{\theta,\kappa}^{\infty}\sum_{\phi\in\theta\cdot\kappa} \frac{h^{(2t+l)}
    (\tr\boldsymbol{\Omega})\Delta_{\phi}^{\theta,\kappa}
    C_{\phi}\left(\mathbf{D}^{2}\right)C_{\phi}^{\theta,\kappa}
    \left(\boldsymbol{\Sigma}^{-1},\boldsymbol{\Omega}\boldsymbol{\Sigma}^{-1}\right)}{t!l!
    \left(\frac{1}{2}K\right)_{\kappa}C_{\phi}\left(\mathbf{I}_{N-1}\right)},\\
    &=&\sum_{\theta,\kappa}^{\infty}\sum_{\phi\in\theta\cdot\kappa}
    \frac{h^{(2t+l)}(\tr\boldsymbol{\Omega})\left(\Delta_{\phi}^{\theta,\kappa}\right)^{2}
    C_{\phi}\left(\mathbf{D}^{2}\right)C_{\kappa} \left(\boldsymbol{\Omega}\right)}{t!l!
    \left(\sigma^{2}\right)^{2t+l}\left(\frac{1}{2}K\right)_{\kappa}
    C_{\kappa}\left(\mathbf{I}_{N-1}\right)}.
\end{eqnarray*}
Note that $\sum_{\phi\in\theta\cdot\kappa}
\left(\Delta_{\phi}^{\theta,\kappa}\right)^{2}C_{\phi}\left(\mathbf{D}^{2}\right)=C_{\kappa}\left(\mathbf{D}^{2}\right)C_{\theta}\left(\mathbf{D}^{2}\right)$,
see \citet{d:80}, eq. (5.10). Thus
\begin{eqnarray*}
    J&=&\sum_{\theta,\kappa}^{\infty} \frac{h^{(2t+l)}(\tr\boldsymbol{\Omega})
    C_{\kappa}\left(\mathbf{D}^{2}\right)C_{\theta}\left(\mathbf{D}^{2}\right)C_{\kappa}
    \left(\boldsymbol{\Omega}\right)}{t!l!\left(\sigma^{2}\right)^{2t+l}
    \left(\frac{1}{2}K\right)_{\kappa}C_{\kappa}\left(\mathbf{I}_{N-1}\right)},\\
    &=& \sum_{\kappa}^{\infty}\sum_{l=0}^{\infty}
    \frac{h^{(2t+l)}(\tr\boldsymbol{\Omega})
    C_{\kappa}\left(\mathbf{D}^{2}\right)C_{\kappa}
    \left(\boldsymbol{\Omega}\right)}{t!l!\left(\sigma^{2}\right)^{2t+l}
    \left(\frac{1}{2}K\right)_{\kappa}C_{\kappa}\left(\mathbf{I}_{N-1}\right)}
    \sum_{\theta}C_{\theta}\left(\mathbf{D}^{2}\right),\\
    &=&\sum_{\kappa}^{\infty}\sum_{l=0}^{\infty} \frac{h^{(2t+l)}(\tr\boldsymbol{\Omega})
    C_{\kappa}\left(\mathbf{D}^{2}\right)C_{\kappa} \left(\boldsymbol{\Omega}\right)
    }{t!l!\left(\sigma^{2}\right)^{2t+l}\left(\frac{1}{2}K\right)_{\kappa}
    C_{\kappa}\left(\mathbf{I}_{N-1}\right)} \left(\tr \mathbf{D}^{2}\right)^{l}.
\end{eqnarray*}
Now, observe that
$h(v)=\sum_{l=0}^{\infty}\frac{h^{(l)}(a)}{l!}(v-a)^{l}$, with
$a=\tr \boldsymbol{\Omega}$, $v=\tr\boldsymbol{\Omega}+\frac{1}{\sigma^{2}}\tr \mathbf{D}^{2}$, and
$h(v)=h^{(2t)}(v)$. Thus
\begin{eqnarray*}
    J &=& \sum_{\kappa}^{\infty} \frac{h^{(2t)}(\tr \boldsymbol{\Omega} +
    \frac{1}{\sigma^{2}}\tr \mathbf{D}^{2})
    C_{\kappa}\left(\mathbf{D}^{2}\right)C_{\kappa} \left(\boldsymbol{\Omega}\right)
    }{t!\left(\sigma^{2}\right)^{2t} \left(\frac{1}{2}K\right)_{\kappa}C_{\kappa}
    \left(\mathbf{I}_{N-1}\right)},\\
    &=& \sum_{\kappa}^{\infty} \frac{h^{(2t)}(\tr \boldsymbol{\Omega} +
    \frac{1}{\sigma^{2}} \tr \mathbf{D}^{2}) C_{\kappa}\left(\frac{1}{\sigma^{2}}
    \mathbf{D}^{2}\right)C_{\kappa} \left(\boldsymbol{\Omega}\right)
    }{t!\left(\frac{1}{2}K\right)_{\kappa}C_{\kappa}\left(\mathbf{I}_{N-1}\right)}. \qed
\end{eqnarray*}

Now let be $\mathbf{W} = \mathbf{D}/r$, $r= \|\mathbf{D}\| = \|\mathbf{V}' \mathbf{D}
\mathbf{\mathbf{H}}\| = \|\mathbf{Y}\|$ and noting that if $\mathbf{D} = \diag(D_{1},
\dots, D_{n})$ we define $\vecp(\mathbf{D}) = (D_{1}, \dots, D_{n})$, then
$$
\vecp (\mathbf{D})=\left(%
\begin{array}{c}
  D_{1} \\
  \vdots \\
  D_{n} \\
\end{array}%
\right) \mbox{, implies that }
\vecp (\mathbf{W})=\left(%
\begin{array}{c}
  D_{1}/r \\
  \vdots \\
  D_{n}/r \\
\end{array}%
\right)=\vecp \frac{\mathbf{D}}{r},
$$
thus
\begin{eqnarray*}
(d\mathbf{W}(\mathbf{u}))&=&r^{m}\prod_{i=1}^{m}\sin^{m-i}\theta_{i}(d\mathbf{u})\wedge dr\\
&=&r^{m}J(\mathbf{u})(d\mathbf{u})\wedge dr,
\end{eqnarray*}
with $m=n-1$, $\mathbf{u}=\left(\theta_{1},\ldots,\theta_{m}\right)'$.

The shape density under Le and Kendall's approach it is known as disk density.

\begin{thm}\label{th:diskdensity}
The disk density is given by
\begin{eqnarray}
f_{\mathbf{W}}(\mathbf{W})&=&\frac{\pi^{\frac{nK}{2}}\displaystyle
\prod_{i=1}^{n}l_{i}^{*N-1+K-2n}\prod_{i<j}\left(l_{i}^{*2}-l_{j}^{*2}\right)J(\mathbf{u})}
{\Gamma_{n}\left[\frac{1}{2} K\right]| \boldsymbol{\Sigma}|^{\frac{K}{2}}}\nonumber\\
&& \times \sum_{t=0}^{\infty}\sum_{\kappa}
\frac{1}{t!\left(\frac{1}{2}K\right)_{\kappa}} \int_{\mathbf{V} \in
\mathcal{V}_{n,N-1}} \frac{C_{\kappa}\left( \mathbf{\Omega} \mathbf{\Sigma}^{-1}
\mathbf{V}'\mathbf{W}^{2}\mathbf{V}\right)} {\displaystyle
\sum_{\theta}C_{\theta}\left(\mathbf{\Sigma}^{-1}\mathbf{V}'\mathbf{W}^{2}\mathbf{V}\right)}
(\mathbf{V}d\mathbf{V}')\nonumber\\
&& \label{Teo:generaldiskdensity}\times
\int_{0}^{\infty}s^{n(N+K-n-1)+2t-1}h^{(2t)}\left(s^{2} + \tr
\mathbf{\Omega}\right)(ds),
\end{eqnarray}
where the number of landmarks $N$ are selected in such way that
$\frac{n(N+K-n-1)}{2}+t$ is a positive integer, then
$\theta=(l_{i},\ldots,l_{\alpha})$ is a partition of $\frac{n(N+K-n-1)}{2}+t$, and
$\sum_{i=1}^{\alpha}l_{i}=\frac{n(N+K-n-1)}{2}+t$
\end{thm}
\textit{Proof.} From Theorem \ref{th:jointVD}
\begin{eqnarray*}
    dF_{\mathbf{V},\mathbf{D}}(\mathbf{V},\mathbf{D})&=&\frac{\pi^{\frac{nK}{2}}
    |\mathbf{D}|^{N-1+K-2n}\displaystyle\prod_{i<j}\left(D_{i}^{2}-D_{j}^{2}\right)}
    {\Gamma_{n}\left[\frac{1}{2} K\right]|\boldsymbol{\Sigma}|^{\frac{K}{2}}}(d\mathbf{D})
    (\mathbf{V}d\mathbf{V}')\\
    &&\times \sum_{t=0}^{\infty}\sum_{\kappa}\frac{h^{(2t)} \left[\tr
    \left(\boldsymbol{\Sigma}^{-1}\mathbf{V}'\mathbf{D}^{2}\mathbf{V}\right)\right]
    C_{\kappa}\left(\boldsymbol{\Omega}\boldsymbol{\Sigma}^{-1}\mathbf{V}'
    \mathbf{D}^{2}\mathbf{V}\right)}{t!\left(\frac{1}{2}K\right)_{\kappa}}.
\end{eqnarray*}
Let $\mathbf{W}=\diag (l_{1}^{*},\ldots,l_{n}^{*})$, $l_{i}^{*}=\frac{D_{i}}{r}$,
$r=\|\mathbf{D}\|=\|\mathbf{Y}\|$, then

\begin{eqnarray*}
    dF_{\mathbf{V},\mathbf{W}}(\mathbf{V},\mathbf{W})&=&\frac{\pi^{\frac{nK}{2}}
    |r\mathbf{W}|^{N-1+K-2n}r^{m}J(\mathbf{u}) \displaystyle
    \prod_{i<j}r^{2}\left(l_{i}^{*2}-l_{j}^{*2} \right)}{\Gamma_{n}\left[\frac{1}{2}
    K\right]|\boldsymbol{\Sigma}|^{\frac{K}{2}}}(d\mathbf{W})(\mathbf{V}d\mathbf{V}')\\
    &&\times \sum_{t=0}^{\infty}\sum_{\kappa}\frac{h^{(2t)} \left[r^{2}\tr
    \left(\boldsymbol{\Sigma}^{-1}\mathbf{V}' \mathbf{W}^{2}\mathbf{V}\right)\right]
    C_{\kappa}\left(r^{2}\boldsymbol{\Omega}\boldsymbol{\Sigma}^{-1}\mathbf{V}'
    \mathbf{W}^{2}\mathbf{V}\right)}{t!\left(\frac{1}{2}K\right)_{\kappa}}.
\end{eqnarray*}
Note that
\begin{enumerate}
    \item $|r\mathbf{W}|^{N-1+K-2n}=r^{n(N-1+K-2n)}$.
    \item $\displaystyle\prod_{i<j}r^{2}\left(l_{i}^{*2}-l_{j}^{*2}\right) =
    \left(r^{2}\right)^{\frac{n(n-1)}{2}}\displaystyle\prod_{i<j}
    \left(l_{i}^{*2}-l_{j}^{*2}\right)$.
    \item
    $C_{\kappa}\left(r^{2}\boldsymbol{\Omega}\boldsymbol{\Sigma}^{-1}\mathbf{V}'
    \mathbf{W}^{2}\mathbf{V}\right)=r^{2t}C_{\kappa}\left(\boldsymbol{\Omega}
    \boldsymbol{\Sigma}^{-1}\mathbf{V}'\mathbf{W}^{2}\mathbf{V}\right)$.
\end{enumerate}
Collection powers of $r$ and defining $r = \frac{s}{\left(\tr
\boldsymbol{\Sigma}^{-1} \mathbf{V}' \mathbf{W}^{2}\mathbf{V}
\right)^{\frac{1}{2}}}$, with $dr = \frac{ds}{\left(\tr\boldsymbol{\Sigma}^{-1}
\mathbf{V}'\mathbf{W}^{2}\mathbf{V}\right)^{\frac{1}{2}}}$ and
\begin{footnotesize}
\begin{eqnarray*}
    &&\int_{0}^{\infty}r^{n(N+K-n-1)+2t-1}h^{(2t)}\left(r^{2} \tr
    \boldsymbol{\Sigma}^{-1}\mathbf{V}' \mathbf{W}^{2}\mathbf{V}\right)dr\\
    &&\quad= \int_{0}^{\infty} \left(\frac{s}{\left(\tr\boldsymbol{\Sigma}^{-1}
    \mathbf{V}'\mathbf{W}^{2}\mathbf{V}\right)^{\frac{1}{2}}}\right)^{n(N+K-n-1)+2t-1}
    h^{(2t)}\left(s^{2}\right)\frac{ds}{\left(\tr\boldsymbol{\Sigma}^{-1}\mathbf{V}'
    \mathbf{W}^{2}\mathbf{V}\right)^{\frac{1}{2}}}\\
    && \quad=\left(\tr\boldsymbol{\Sigma}^{-1}\mathbf{V}'\mathbf{W}^{2}
    \mathbf{V}\right)^{-\frac{n(N+K-n-1)}{2}+t} \int_{0}^{\infty}s^{n(N+K-n-1)+2t-1}
    h^{(2t)}\left(s^{2}\right)ds.
\end{eqnarray*}
\end{footnotesize}
Thus the marginal density of $dF_{\mathbf{W}}(\mathbf{W})$ is given by
\begin{eqnarray*}
    &=&\frac{\pi^{\frac{nK}{2}}|\mathbf{W}|^{N-1+K-2n}J(\mathbf{u})
    \displaystyle\prod_{i<j}\left(l_{i}^{*2}-l_{j}^{*2}\right)} {\Gamma_{n}
    \left[\frac{K}{2}\right]|\boldsymbol{\Sigma}|^{\frac{K}{2}}}
    \sum_{t=0}^{\infty}\sum_{\kappa}\frac{1}{t!\left(\frac{1}{2}K\right)_{\kappa}}\\
    && \times \int_{\mathbf{V}\in \mathcal{V}_{n,N-1}}
    C_{\kappa}\left(\boldsymbol{\Omega}\boldsymbol{\Sigma}^{-1} \mathbf{V}'
    \mathbf{W}^{2} \mathbf{V}\right)\left(\tr\boldsymbol{\Sigma}^{-1}\mathbf{V}'
    \mathbf{W}^{2} \mathbf{V}\right)^{-\frac{n(N+K-n-1)}{2}+t}(\mathbf{V}d\mathbf{V}')\\
    &&\times \int_{0}^{\infty}s^{n(N+K-n-1)+2t-1} h^{(2t)}\left(s^{2}\right)ds.
\end{eqnarray*}
Now, let be
\begin{eqnarray*}
    J &=& \int_{\mathbf{V}\in V_{n,N-1}}\frac{C_{\kappa} \left(\boldsymbol{\Omega}
    \boldsymbol{\Sigma}^{-1} \mathbf{V}'\mathbf{W}^{2} \mathbf{V}\right)
    (\mathbf{V}d\mathbf{V}')}{\left(\tr \boldsymbol{\Sigma}^{-1}\mathbf{V}'
    \mathbf{W}^{2} \mathbf{V} \right)^{\frac{n(N+K-n-1)}{2}+t}},\\
    &=& \int_{\mathbf{V} \in V_{n,N-1}} \frac{C_{\kappa}\left( \boldsymbol{\Omega}
    \boldsymbol{\Sigma}^{-1}\mathbf{V}' \mathbf{W}^{2}\mathbf{V}\right)}{\displaystyle
    \sum_{\theta}C_{\theta}\left(\boldsymbol{\Sigma}^{-1} \mathbf{V}' \mathbf{W}^{2}
    \mathbf{V}\right)}(\mathbf{V}d\mathbf{V}'),
\end{eqnarray*}
where the number of landmarks $N$ are selected in such way that $\frac{n(N+K-n-1)}{2}
+ t$ is a positive integer, then $\theta=(l_{i},\ldots,l_{\alpha})$ is a partition of
$\frac{n(N+K-n-1)}{2}+t$, and $\sum_{i=1}^{\alpha}l_{i}=\frac{n(N+K-n-1)}{2}+t$. Then
we obtain the desired result. \qed

The isotropic case of the disk distribution follows

\begin{thm}\label{th:diskIsotropic}
The isotropic disk density is given by
\begin{eqnarray}
    f_{\mathbf{W}}(\mathbf{W})&=&\frac{2^{n}\pi^{\frac{n(N+K-1)}{2}}
    \displaystyle\prod_{i=1}^{n}l_{i}^{*\,N-1+K-2n}\prod_{i<j}
    \left(l_{i}^{*2}-l_{j}^{*2}\right)J(\mathbf{u})}
    {\Gamma_{n}\left[\frac{1}{2} K\right]\Gamma_{n}
    \left[\frac{N-1}{2}\right]\left(\sigma^{2}\right)^{\frac{(N-1)K}{2}}}\nonumber\\
    && \times \sum_{\kappa}^{\infty} \frac{C_{\kappa}
    \left(\displaystyle\frac{1}{\sigma^{2}}\mathbf{W}^{2}\right)
    C_{\kappa}\left(\boldsymbol{\Omega}\right)} {t!\left(\frac{1}{2}K\right)_{\kappa}
    C_{\kappa}\left(\mathbf{I}_{N-1}\right)}\nonumber\\
    && \label{Coro:Isodickdensity} \times \int_{0}^{\infty}r^{n(N+K-n-1)+2t-1}h^{(2t)}
    \left(\tr\boldsymbol{\Omega}+\frac{r^{2}}{\sigma^{2}}\right)dr.
\end{eqnarray}
\end{thm}
 \textit{Proof.} The result is obtained from
(\ref{Teo:generaldiskdensity}) taking $\mathbf{\Sigma} = \sigma^{2}\mathbf{I}$ and
observing that $|\mathbf{\Sigma}| = |\sigma^{2}\mathbf{I}| = (\sigma^{2})^{(N-1)}$;
$$
  \tr\boldsymbol{\Sigma}^{-1}\mathbf{V}'\mathbf{W}^{2}\mathbf{V} = \frac{1}{\sigma^{2}}
  \tr\mathbf{V}'\mathbf{W}^{2} \mathbf{V} = \frac{1}{\sigma^{2}}\tr\mathbf{W}^{2}
  \mathbf{V}\mathbf{V}' =\frac{1}{\sigma^{2}}\tr\mathbf{W}^{2} =
  \frac{1}{\sigma^{2}},
$$
recalling that $||\mathbf{W}|| = 1$. Hence
\begin{eqnarray*}
    J &=&\int_{\mathbf{V}\in V_{n,N-1}}\frac{C_{\kappa}
    \left(\boldsymbol{\Omega}\boldsymbol{\Sigma}^{-1}\mathbf{V}'\mathbf{W}^{2}\mathbf{V}\right)
    (\mathbf{V}d\mathbf{V}')}{\left(\tr\boldsymbol{\Sigma}^{-1}\mathbf{V}' \mathbf{W}^{2}
    \mathbf{V} \right)^{\frac{n(N+K-n-1)}{2}+t}}\\
    &=& \left(\sigma^{2}\right)^{n(N+K-n-1)+t}\int_{\mathbf{V} \in V_{n,N-1}}
    C_{\kappa}\left(\frac{1}{\sigma^{2}}\boldsymbol{\Omega}\mathbf{V}'
    \mathbf{W}^{2}\mathbf{V}\right)(\mathbf{V}d\mathbf{V}')\\
    &=& \frac{\left(\sigma^{2}\right)^{n(N+K-n-1)+t} 2^{n} \pi^{n(N-1)/2}
    }{\Gamma_{n}\left[ \frac{1}{2}(N-1)\right]} \frac{C_{\kappa}\left(\displaystyle
    \frac{1}{\sigma^{2}} \mathbf{\Omega}\right) C_{\kappa}\left(\mathbf{W}^{2}\right)}
    {C_{\kappa}\left(\mathbf{I}_{N-1}\right)}.
\end{eqnarray*}
Finally, (\ref{Coro:Isodickdensity}) is obtained making the change of variable $s =
r/\sigma$ with $ds = dr/\sigma$ in (\ref{Teo:generaldiskdensity}) and observing that
$C_{\kappa}\left(\displaystyle\frac{1}{\sigma^{2}} \mathbf{\Omega}\right)
C_{\kappa}\left(\mathbf{W}^{2}\right) = C_{\kappa}\left( \mathbf{\Omega}\right)
C_{\kappa}\left(\displaystyle\frac{1}{\sigma^{2}}\mathbf{W}^{2}\right)$.

Alternatively, let $l_{i}^{*}= D_{i}/r$, $\tr \mathbf{D}^{2}=\sum_{i=1}^{n}
D_{i}^{2}=r^{2}$ in (\ref{eq:coneIsotropic}), therefore
\begin{eqnarray*}
    f_{\mathbf{W}}(\mathbf{W})&=&\frac{2^{n}\pi^{\frac{n(N+K-1)}{2}}r^{m}
    J(\mathbf{u})\prod_{i=1}^{n}\left(r\,l_{i}^{*}\right)^{N-1+K-2n}
    \displaystyle\prod_{i<j}r^{2}\left(l_{i}^{*2}-l_{j}^{*2}\right)}
    {\Gamma_{n}\left[\frac{1}{2} K\right]\Gamma_{n}\left[\frac{1}{2}
    (N-1)\right]\left(\sigma^{2}\right)^{\frac{(N-1)K}{2}}}\\
    && \times \sum_{t=0}^{\infty}\sum_{\kappa} \frac{h^{(2t)} \left(\tr
    \boldsymbol{\Omega} + \frac{r^{2}}{\sigma^{2}}\right)
    C_{\kappa}\left(\frac{r^{2}}{\sigma^{2}} \mathbf{W}^{2}\right)
    C_{\kappa}\left(\boldsymbol{\Omega}\right)} {t!\left(\frac{1}{2}K\right)_{\kappa}
    C_{\kappa}\left(\mathbf{I}_{N-1}\right)}.
\end{eqnarray*}
Observe that
\begin{enumerate}
    \item
    $\prod_{i=1}^{n}\left(r\,l_{i}^{*}\right)^{N-1+K-2n}=r^{n(N-1+K-2n)}
    \prod_{i=1}^{n}l_{i}^{*\,N-1+K-2n}$.
    \item $\displaystyle\prod_{i<j}r^{2}\left(l_{i}^{*2}-l_{j}^{*2}\right) =
    r^{n(n-1)}\displaystyle\prod_{i<j}\left(l_{i}^{*2}-l_{j}^{*2}\right)$
    \item
    $C_{\kappa}\left(\frac{r^{2}}{\sigma^{2}}\mathbf{W}^{2}\right) =
    r^{2t}C_{\kappa}\left(\frac{1}{\sigma^{2}}\mathbf{W}^{2}\right)$.
\end{enumerate}
Collecting powers of $r$ we have
\begin{equation*}
    \int_{0}^{\infty}r^{n(N+K-n-1)+2t-1}h^{(2t)}\left(\tr\boldsymbol{\Omega} +
    \frac{r^{2}}{\sigma^{2}}\right)dr. \qed
\end{equation*}

\section{Central Case}\label{sec:Centralconedisk}

We obtain the central cases of the cone and the disk densities.
\begin{cor}\label{cor:CentralCone}
The central cone density is given by
\begin{small}
\begin{equation}\label{eq:CentralCone}
    f_{\mathbf{D}}(\mathbf{D})=\frac{2^{n}|\mathbf{D}|^{N-1+K-2n}\displaystyle
    \prod_{i<j}\left(D_{i}^{2}-D_{j}^{2}\right)} {\pi^{-\frac{n(N+K-1)}{2}}\Gamma_{n}
    \left[\frac{1}{2} K\right]\Gamma_{n}\left[\frac{1}{2} (N-1)\right]
    |\boldsymbol{\Sigma}|^{\frac{K}{2}}}\sum_{l=0}^{\infty}
    \sum_{\theta}\frac{h^{(l)}(0)C_{\theta}\left(\boldsymbol{\Sigma}^{-1}\right)
    C_{\theta}\left(\mathbf{D}^{2}\right)}{l!C_{\theta}\left(\mathbf{I}_{N-1}\right)}.
\end{equation}
\end{small}
\end{cor}
\textit{Proof.} Start with
$$
  dF_{\mathbf{Y}}(\mathbf{Y})=\frac{1}{|\boldsymbol{\Sigma}|^{\frac{K}{2}}} h
  \left[\tr\boldsymbol{\Sigma}^{-1}\mathbf{Y}\mathbf{Y}'\right](d\mathbf{Y}).
$$
The joint density of $\mathbf{V}, \mathbf{D}, \mathbf{\mathbf{H}}$ is
$$
  dF_{\mathbf{V},\mathbf{D},\mathbf{\mathbf{H}}}(\mathbf{V}, \mathbf{D},
  \mathbf{\mathbf{H}}) = \frac{\displaystyle\prod_{i<j}
  \left(D_{i}^{2}-D_{j}^{2}\right)(d\mathbf{D})} {2^{n}|\mathbf{D}|^{-(N-1+K-2n)}
  |\boldsymbol{\Sigma}|^{\frac{K}{2}}} h \left[\tr\boldsymbol{\Sigma}^{-1} \mathbf{V}'
  \mathbf{D}^{2} \mathbf{V}\right]
  (\mathbf{\mathbf{H}}d\mathbf{\mathbf{H}}')(\mathbf{V}d\mathbf{V}').
$$
Recalling that
$$
  \int_{\mathbf{\mathbf{H}} \in \mathcal{V}_{n,K}}(\mathbf{\mathbf{H}}'
  d\mathbf{\mathbf{H}})=\frac{2^{n}\pi^{\frac{nK}{2}}}{\Gamma_{n}\left[\frac{1}{2}
  K\right]}.
$$
we have that
$$
  dF_{\mathbf{V},\mathbf{D}}(\mathbf{V},\mathbf{D}) = \frac{\pi^{\frac{nK}{2}} |\mathbf{D}|^{N-1+K-2n}
  \displaystyle\prod_{i<j}\left(D_{i}^{2}-D_{j}^{2}\right)(d\mathbf{D})}
  {\Gamma_{n}\left[\frac{1}{2}K\right]|\boldsymbol{\Sigma}|^{\frac{K}{2}}} h \left[
  \tr\boldsymbol{\Sigma}^{-1}\mathbf{V}'\mathbf{D}^{2}\mathbf{V}\right](\mathbf{V}d\mathbf{V}').
$$
Then the required density is given by%
$$
  f_{\mathbf{D}}(\mathbf{D})=\frac{|\mathbf{D}|^{N-1+K-2n} \displaystyle \prod_{i<j}
  \left(D_{i}^{2}-D_{j}^{2}\right)}{\pi^{-\frac{nK}{2}}\Gamma_{n}\left[\frac{1}{2}
  K\right]|\boldsymbol{\Sigma}|^{\frac{K}{2}}}\int_{\mathbf{V}\in \mathcal{V}_{n,N-1}}
  h \left[\tr\boldsymbol{\Sigma}^{-1} \mathbf{V}' \mathbf{D}^{2} \mathbf{V}\right]
  (\mathbf{V}d\mathbf{V}').
$$
Integrating with respect to $\mathbf{V}$
\begin{eqnarray*}
    &&\displaystyle\int_{\mathcal{V}_{n,N-1}}h\left[\tr\boldsymbol{\Sigma}^{-1}\mathbf{V}'
    \mathbf{D}^{2} \mathbf{V}\right](\mathbf{V}d\mathbf{V}')\\
    && \hspace{3cm}  = \int_{\mathcal{V}_{n,N-1}}\sum_{l=0}^{\infty}
    \frac{h^{(l)}(0)}{l!}\left(\tr\boldsymbol{\Sigma}^{-1}\mathbf{V}'\mathbf{D}^{2}
    \mathbf{V}\right)^{l}(\mathbf{V}d\mathbf{V}')\\
    && \hspace{3cm}=\sum_{l=0}^{\infty}\sum_{\theta}\frac{h^{(l)}(0)}{l!}\int_{
    \mathcal{V}_{n,N-1}}C_{\theta}\left(\boldsymbol{\Sigma}^{-1}\mathbf{V}'
    \mathbf{D}^{2}\mathbf{V}\right)(\mathbf{V}d\mathbf{V}')\\
    &&\hspace{3cm}=\sum_{l=0}^{\infty}\sum_{\theta}\frac{h^{(l)}(0)}{l!}
    \frac{2^{n}\pi^{\frac{n(N-1)}{2}}C_{\theta}\left(\boldsymbol{\Sigma}^{-1}\right)
    C_{\theta}\left(\mathbf{D}^{2}\right)} {\Gamma_{n}\left[\frac{1}{2}
    (N-1)\right]C_{\theta}\left(\mathbf{I}_{N-1}\right)}\\
    &&\hspace{3cm}=\frac{2^{n}\pi^{\frac{n(N-1)}{2}}}{\Gamma_{n}\left[\frac{1}{2}
    (N-1)\right]} \sum_{l=0}^{\infty}\sum_{\theta}\frac{h^{(l)}(0)}{l!}
    \frac{C_{\theta}\left(\boldsymbol{\Sigma}^{-1}\right)C_{\theta}\left(\mathbf{D}^{2}\right)}
    {C_{\theta}\left(\mathbf{I}_{N-1}\right)}.
\end{eqnarray*}
Then,
\begin{eqnarray*}
    f_{\mathbf{D}}(\mathbf{D})&=& \frac{2^{n} \pi^{\frac{n(N+K-1)}{2}}
    |\mathbf{D}|^{N-1+K-2n}\displaystyle \prod_{i<j}\left(D_{i}^{2}-D_{j}^{2}\right)} {
    \Gamma_{n}\left[\frac{1}{2} K\right]\Gamma_{n}\left[\frac{1}{2}
    (N-1)\right]|\boldsymbol{\Sigma}|^{\frac{K}{2}}}\\
    && \hspace{4cm}\times \sum_{l=0}^{\infty}\sum_{\theta}\frac{h^{(l)}(0)}{l!}
    \frac{C_{\theta}\left(\boldsymbol{\Sigma}^{-1}\right)C_{\theta}\left(\mathbf{D}^{2}\right)}
    {C_{\theta}\left(\mathbf{I}_{N-1}\right)}.
\end{eqnarray*}

Alternatively, from Theorem \ref{th:conedensity}, if we take
\begin{enumerate}
  \item $h^{(2t+l)}\left(\tr\boldsymbol{\Omega}\right)= h^{(l)}(0),$
  \item $\Delta_{\phi}^{\theta,\kappa}=\Delta_{\phi}^{\theta,\theta}=
        \displaystyle\frac{C_{\theta}(\mathbf{I})}{C_{\theta}(\mathbf{I})}=1,$
  \item $C_{\phi}(\mathbf{D}^{2}) = C_{\theta}(\mathbf{D}^{2}),
        \quad C_{\phi}(\mathbf{I})=C_{\theta}(\mathbf{I}),$
  \item $C_{\phi}^{\theta,\kappa}\left(\boldsymbol{\Sigma}^{-1},
        \boldsymbol{\Omega}\boldsymbol{\Sigma}^{-1}\right)=
        C_{\phi}^{\theta,\kappa}\left(\boldsymbol{\Sigma}^{-1},0\right)=
        C_{\theta}\left(\boldsymbol{\Sigma}^{-1}\right),$
\end{enumerate}
the required result follows. \qed

And finally, note that the central disk density is invariant under the elliptical
distributions.
\begin{cor}\label{cor:Centraldisk}
The central disk density is given by
\begin{eqnarray}
    f_{\mathbf{W}}(\mathbf{W})&=&\frac{\Gamma[\frac{1}{2}(n(N+K-n-1))]\displaystyle
    \prod_{i=1}^{n} l_{i}^{*N-1+K-2n} \prod_{i<j}\left(l_{i}^{*2}-l_{j}^{*2}\right)
    J(\mathbf{u})}{2 \ \pi^{\frac{n(N-n-1)}{2}} \Gamma_{n} \left[\frac{1}{2} K\right]|
    \boldsymbol{\Sigma}|^{\frac{K}{2}}}\nonumber\\
    && \label{centralnoisodisk}\hspace{4cm}\times \int_{\mathbf{V} \in
    \mathcal{V}_{n,N-1}} \frac{(\mathbf{V}d\mathbf{V}')} {\displaystyle
    \sum_{\theta}C_{\theta}\left(\mathbf{\Sigma}^{-1}\mathbf{V}'\mathbf{W}^{2}
    \mathbf{V}\right)}.
\end{eqnarray}
where the number of landmarks $N$ are selected in such way that
$\frac{n(N+K-n-1)}{2}$ is a positive integer, then $\theta=(l_{i},\ldots,l_{\alpha})$
is a partition of $\frac{n(N+K-n-1)}{2}$, and
$\sum_{i=1}^{\alpha}l_{i}=\frac{n(N+K-n-1)}{2}$ where
$\theta=(l_{i},\ldots,l_{\alpha})$ is a partition of the positive integer
$n(N+K-n-1)/2$, $\displaystyle\sum_{i=1}^{\alpha}l_{i}=n(N+K-n-1)/2$.
\end{cor}
\textit{Proof.} From Theorem \ref{th:diskdensity}, taking $t=0$,
$\boldsymbol{\Omega}=\mathbf{0}$, $h^{(0)}(\cdot)\equiv h(\cdot)$ and recalling that
$$
  \int_{0}^{\infty}s^{n(N+K-n-1)-1}h\left(s^{2}\right)(ds) =
  \frac{\Gamma[\frac{1}{2}(n(N+K-n-1))]}{\pi^{\frac{n(N+K-n-1)}{2}}},
$$
the result is obtained.\qed

\section{Example: Mouse Vertebra}\label{sub:mouse}

As the reader can check the general cone and disk densities are given in terms of
invariant polynomials, so at this time no inference can be performed, except if the
series are truncated in the first few terms.  However, there is a way to work with an
exact density, it is, when we assume an isotropic model.

Consider the isotropic elliptical disk density of theorem \ref{th:diskIsotropic}
\begin{eqnarray*}
    f_{\mathbf{W}}(\mathbf{W})&=&\frac{2^{n}\pi^{\frac{n(N+K-1)}{2}}\displaystyle
    \prod_{i=1}^{n}l_{i}^{*\,N-1+K-2n}\prod_{i<j}\left(l_{i}^{*2}-l_{j}^{*2}\right)
    J(\mathbf{u})} {\Gamma_{n}\left[\frac{1}{2} K\right]\Gamma_{n}\left[\frac{1}{2}
    (N-1)\right]\left(\sigma^{2}\right)^{\frac{(N-1)K}{2}}}\\
    &&\times \sum_{\kappa}^{\infty} \frac{C_{\kappa}
    \left(\frac{1}{\sigma^{2}}\mathbf{W}^{2}\right)
    C_{\kappa}\left(\boldsymbol{\Omega}\right)} {t!\left(\frac{1}{2}K\right)_{\kappa }
    C_{\kappa}\left(\mathbf{I}_{N-1}\right)}\\
    &&\times \int_{0}^{\infty}r^{n(N+K-n-1)+2t-1}h^{(2t)}
    \left(\tr\boldsymbol{\Omega}+\frac{r^{2}}{\sigma^{2}}\right)dr,
\end{eqnarray*}
and the generator for the subfamily Kotz
$$
  h(y) = \frac{R^{T-1+\frac{K(N-1)}{2}}\Gamma\left(\frac{K(N-1)}{2}\right)}{\pi^{K(N-1)/2}
  \Gamma\left(T-1+\frac{K(N-1)}{2}\right)}y^{T-1}\e^{-Ry},
$$
with derivative
$$
 \frac{d^{k}}{dy^{k}}y^{T-1}\e^{-Ry} = (-R)^{k}y^{T-1}\e^{-Ry} \left\{1+\sum_{m=1}^{k}
 \binom{k}{m}\left[\prod_{i=0}^{m-1}(T-1-i)\right](-Ry)^{-m}\right\},
$$
see \citet{Caro2009} for other families (Pearson VII, Bessel, general Kotz, Jensen
Logistic) and their derivatives.

Now, we contrast three models next, the classical Gaussian ($T=1$, $R=\frac{1}{2}$)
and two non normal ($T=2$, $R=\frac{1}{2}$ and $T=3$, $R=\frac{1}{2}$) via the
modified BIC criterion, they will be applied to the data of two groups (small and
large) of mouse vertebra, and experiment very detailed in \citet{DM98}.

The isotropic Gaussian disk density is given by

\begin{eqnarray*}
    &&f_{\mathbf{W}}(\mathbf{W})=\frac{2^{\frac{1}{2}(-2-M+n+Kn-n^{2}+nN)} \displaystyle
    \prod_{i=1}^{n}l_{i}^{*\,N-1+K-2n}\prod_{i<j}\left(l_{i}^{*2}-l_{j}^{*2}\right)
    J(\mathbf{u})} {\pi^{\frac{1}{2}(M+n-Kn-nN)}\sigma^{M-n(-1+K-n+N)}
    \Gamma_{n}\left[\frac{1}{2} K\right]\Gamma_{n}\left[\frac{1}{2} (N-1)\right]}\\
    &&\times \etr\left(-\frac{\boldsymbol{\mu}'\boldsymbol{\mu}}{2\sigma^{2}}\right)
    \sum_{t=0}^{\infty}\frac{\Gamma\left[\frac{1}{2}(n(-1+K-n+N))+t\right]}{t!}
    \sum_{\kappa}\frac{C_{\kappa}\left(\mathbf{W}^{2}\right)C_{\kappa}
    \left(\frac{\boldsymbol{\mu}'\boldsymbol{\mu}}{2\sigma^{2}}\right)}
    {\left(\frac{1}{2}K\right)_{\kappa}C_{\kappa}\left(\mathbf{I}_{N-1}\right)},
\end{eqnarray*}
where $M=K(N-1)$ and $n=\min\{ (N-1),K\}$ .

The Kotz disk density when $T=2$ and $R=\frac{1}{2}$ follows after some tedious
simplification, and it is given by

\begin{eqnarray*}
    &&f_{\mathbf{W}}(\mathbf{W})=\frac{2^{\frac{1}{2}(-M+n+Kn-n^{2}+nN)} \displaystyle
    \prod_{i=1}^{n}l_{i}^{*\,N-1+K-2n}\prod_{i<j}\left(l_{i}^{*2}-l_{j}^{*2}\right)J(\mathbf{u})}
    {\pi^{\frac{1}{2}(M+n-Kn-nN)}\sigma^{M-n(-1+K-n+N)}\Gamma_{n}\left[\frac{1}{2}
    K\right]\Gamma_{n}\left[\frac{1}{2} (N-1)\right]M}\\
    && \times \etr\left(-\frac{\boldsymbol{\mu}'\boldsymbol{\mu}}{2\sigma^{2}}\right)
    \sum_{t=0}^{\infty}\frac{1}{t!}\left\{\left(\tr\left(\frac{\boldsymbol{\mu}'
    \boldsymbol{\mu}}{2\sigma^{2}} \right)- 2t\right)
    \Gamma\left[\frac{1}{2}(n(-1+K-n+N))+t\right] \right.\\&&
    \left.+\Gamma\left[\frac{1}{2}(n(-1+K-n+N))+t+1\right]\right\}
    \sum_{\kappa}\frac{C_{\kappa}\left(\mathbf{W}^{2}\right)C_{\kappa}
    \left(\frac{\boldsymbol{\mu}'\boldsymbol{\mu}}{2\sigma^{2}}\right)}
    {\left(\frac{1}{2}K\right)_{\kappa}C_{\kappa}\left(\mathbf{I}_{N-1}\right)}.
\end{eqnarray*}
Finally, the corresponding density for the Kotz model $T=3$, is obtained as:
\begin{eqnarray*}
    &&f_{\mathbf{W}}(\mathbf{W})=\frac{2^{\frac{1}{2}(-2-M+n+Kn-n^{2}+nN)} \displaystyle
    \prod_{i=1}^{n}l_{i}^{*\,N-1+K-2n}\prod_{i<j}\left(l_{i}^{*2}-l_{j}^{*2}\right)}
    {\pi^{\frac{1}{2}(M+n-Kn-nN)}\sigma^{M-n(-1+K-n+N)}\Gamma_{n}\left[\frac{1}{2}
    K\right]\Gamma_{n}\left[\frac{1}{2} (N-1)\right]}\\
    && \times \frac{J(\mathbf{u})}{M(M+2)} \etr \left(-\frac{\boldsymbol{\mu}'
    \boldsymbol{\mu}}{2\sigma^{2}}\right)\\
    && \sum_{t=0}^{\infty}\frac{1}{t!} \left\{\left(-8t+16t^{2}-16t \tr
    \left(\frac{\boldsymbol{\mu}'\boldsymbol{\mu}}{2\sigma^{2}}\right)+ 4\tr^{2}
    \left(\frac{\boldsymbol{\mu}'\boldsymbol{\mu}}{2\sigma^{2}} \right)\right)\right.\\
    && \left. + 4\left(\tr \left(\frac{\boldsymbol{\mu}'\boldsymbol{\mu}}{2\sigma^{2}}
    \right) - 4t \right) \Gamma\left[\frac{1}{2}(n(-1+K-n+N))+t+1\right]\right.\\
    && \left.+4\Gamma\left[\frac{1}{2}(n(-1+K-n+N))+t+2 \right] \right\}
    \sum_{\kappa}\frac{C_{\kappa}\left(\mathbf{W}^{2}\right)C_{\kappa}
    \left(\displaystyle\frac{\boldsymbol{\mu}'\boldsymbol{\mu}}{2\sigma^{2}}\right)}
    {\left(\frac{1}{2}K\right)_{\kappa}C_{\kappa}\left(\mathbf{I}_{N-1}\right)},
\end{eqnarray*}
The likelihood based on exact densities require the computation of the above series,
a carefully comparison with the known hypergeometric of two matrix argument indicates
that these distributions can be obtained by a suitable modification of the algorithms
of \citet{KE06}.

In order to decide which the elliptical model is the best one, different criteria
have been employed for the model selection. We shall consider a modification of the
BIC statistic as discussed in \citet{YY07}, and which was first achieved by
\citet{ri:78} in a coding theory framework. The modified BIC is given by:
$$
  BIC^{*}=-2\mathfrak{L}(\widetilde{\boldsymbol{\boldsymbol{\mu}}},\widetilde{\sigma}^{2},h)
  + n_{p}(\log(n+2)-\log 24),
$$
where
$\mathfrak{L}(\widetilde{\boldsymbol{\boldsymbol{\mu}}},\widetilde{\sigma}^{2},h)$ is
the maximum of the log-likelihood function, $n$ is the sample size and $n_{p}$ is the
number of parameters to be estimated for each particular shape density.

As proposed by \citet{kr:95} and \citet{r:95}, the following selection criteria have
been employed for the model selection.

\begin{table}[ht]  \centering \caption{Grades of evidence
corresponding to values of the $BIC^{*}$ difference.}\label{table2}
\medskip
\renewcommand{\arraystretch}{1}
\begin{center}
  \begin{tabular}{cl}
    \hline
    $BIC^{*}$ difference & Evidence\\
    \hline
    0--2 & Weak\\
    2--6 & Positive \\
    6--10 & Strong\\
    $>$ 10 & Very strong\\
    \hline
  \end{tabular}
\end{center}
\end{table}
The maximum likelihood estimators  for location  parameters associated with the small
and large groups are summarised in the following table:

\begin{table}[ht]  \centering \caption{Maximum likelihood estimators.}\label{table3}
\medskip
\renewcommand{\arraystretch}{1}
\begin{center}
\begin{small}\label{Tab4MR}
\begin{tabular}{|c|c|c|c|c|c|c|c|}
  \hline
  % after \\: \hline or \cline{col1-col2} \cline{col3-col4} ...
  Group& $BIC^{*}$ & $\widetilde{\mu}_{11}$ &$\widetilde{\mu}_{12}$& $\widetilde{\mu}_{21}$
  & $\widetilde{\mu}_{22}$ & $\widetilde{\mu}_{31}$ & $\widetilde{\mu}_{32}$  \\
   & $\build{K:T=2}{G}{K:T=3}$&&&&&&\\
  \hline
  Small & $\build{-536.1662}{-518.3844}{-116.7475}$ & $\build{-16.3534}{2.3931}{-1.1862}$ &
  $\build{-75.7142}{-27.7470}{-19.3374}$ & $\build{40.2639}{34.2670}{41.8450}$
  &$\build{-3.9264}{-16.0397}{-15.4168}$ &$\build{39.6012}{-9.9515}{44.9495}$
  & $\build{-15.9524}{0.4755}{0.3280}$
  \\ \hline
  Large& $\build{-536.9652}{-519.1378}{-115.5510}$ &$\build{-19.3968}{-24.8027}{-27.1256}$
  & $\build{-69.2566}{-63.4131}{-61.2295}$ & $\build{41.9566}{35.8304}{36.5588}$
  & $\build{-1.7093}{-9.0617}{-10.8047}$   & $\build{32.1251}{31.9550}{73.2563}$
  & $\build{-20.2930}{-15.0595}{-6.8034}$  \\
  \hline
\end{tabular}

\medskip

\begin{tabular}{|c|c|c|c|c|c|c|c|}
  \hline
  % after \\: \hline or \cline{col1-col2} \cline{col3-col4} ...
   $\widetilde{\mu}_{41}$ & $\widetilde{\mu}_{42}$ & $\widetilde{\mu}_{51}$
   & $\widetilde{\mu}_{52}$ & $\widetilde{\sigma^{2}}$  \\
  \hline
   $\build{16.0621}{17.2967}{21.6988}$ &$\build{15.5973}{24.2432}{18.9804 }$
   & $\build{-52.0684}{14.8582}{-5.8144}$ &$\build{13.4224}{13.5573}{10.1334}$
   & $\build{106.7591}{33.9149}{26.7961}$
  \\ \hline
  $\build{13.2026}{11.4842}{16.8063}$ &$\build{13.8409}{14.7148}{17.3769}$
  & $\build{-48.3278}{-44.1020}{-46.5415}$ &$\build{16.6046}{21.2210}{23.0966}$ &
  $\build{94.2873}{94.1200}{71.1297}$\\
  \hline
\end{tabular}
\end{small}
\end{center}
\end{table}

According to the modified BIC criterion, the Kotz model with parameters $T=2$,
$R=\frac{1}{2}$ and $s=1$ is the most appropriate, among the three elliptical
densities selected, for modeling the data. There is a very strong difference between
the non normal and the classical Gaussian model in this experiment.

Let $\boldsymbol{\boldsymbol{\mu}}_{1}$ and $\boldsymbol{\boldsymbol{\mu}}_{2}$ be
the mean disk of the small and large groups, respectively. We test equal mean shape
under the best model, and the likelihood ratio (based on
$-2\log\Lambda\approx\chi_{10}^{2}$) for the test
$H_{0}:\boldsymbol{\boldsymbol{\mu}}_{1}=\boldsymbol{\boldsymbol{\mu}}_{2}$ vs
$H_{a}:\boldsymbol{\boldsymbol{\mu}}_{1}\neq\boldsymbol{\boldsymbol{\mu}}_{2}$,
provides the p-value $0.92$, which means that there extremely evidence that the mean
shapes of the two groups are equal.

It is important to note that the general densities derived here apply to any
elliptical model; some classical elliptical densities as  Kotz, Pearson II and VII,
Bessel, Jensen-logistic, can be obtained explicitly and applied, however they demand
the computation of the the $k$-th derivative of the generator elliptical function
$h(\cdot)$. This is not a trivial fact, but for the above mentioned families, the
required formulae are available in \citet{Caro2009}. However, isotropic densities are
more tractable because they are expanded in terms of zonal polynomials, instead of
non isotropic distributions which require some additional conditions on the number of
landmarks in order to obtain a know integral expanded in terms of invariant
polynomials. The general densities expanded in terms of invariant polynomials, seem
non computable at this date for large degrees.

\section*{Acknowledgments}

This research work was supported  by University of Medellin (Medellin, Colombia) and
Universidad Aut\'onoma Agraria Antonio Narro (M\'{e}xico),  joint grant No. 469,
SUMMA group. Also, the first author was partially supported  by CONACYT - M\'exico,
research grant no. \ 138713 and IDI-Spain, Grants No. \ FQM2006-2271 and
MTM2008-05785 and the paper was written during J. A. D\'{\i}az- Garc\'{\i}a's stay as
a visiting professor at the Department of Statistics and O. R. of the University of
Granada, Spain.

\end{document}